\documentclass[12pt]{amsart}
\usepackage{amssymb}

\newcommand{\Bgp}{{\Z^\N}}

\long\def\forget#1\forgotten{}
\newcommand{\issuenumber}{26}
\newcommand{\issuemonth}{December}
\newcommand{\issueyear}{2008}

\newtheorem{issue}{Issue}

\theoremstyle{definition}

\theoremstyle{remark}

\newcommand{\ed}{
\newpage

\section{Unsolved problems from earlier issues}

\begin{issue}
Is $\binom{\Omega}{\Gamma}=\binom{\Omega}{\Tau}$?
\end{issue}

\begin{issue}
Is $\ufin(\cO,\Omega)=\sfin(\Gamma,\Omega)$?
And if not, does $\ufin(\cO,\Gamma)$ imply
$\sfin(\Gamma,\Omega)$?
\end{issue}

\stepcounter{issue}

\begin{issue}
Does $\sone(\Omega,\Tau)$ imply $\ufin(\Gamma,\Gamma)$?
\end{issue}

\begin{issue}
Is $\fp=\fp^*$? (See the definition of $\fp^*$ in that issue.)
\end{issue}

\begin{issue}
Does there exist (in ZFC) an uncountable set satisfying $\sfin(\B,\B)$?
\end{issue}

\stepcounter{issue}

\begin{issue}
Does $X \nin \NON(\M)$ and $Y\nin\mathsf{D}$ imply that
$X\cup Y\nin \COF(\M)$?
\end{issue}

\begin{issue}[CH]
Is $\split(\Lambda,\Lambda)$ preserved under finite unions?
\end{issue}

\begin{issue}
Is $\cov(\M)=\fo$? (See the definition of $\fo$ in that issue.)
\end{issue}

\begin{issue}
Does $\sone(\Gamma,\Gamma)$ always contain an element of cardinality $\fb$?
\end{issue}

\begin{issue}
Could there be a Baire metric space $M$ of weight $\aleph_1$ and a partition
$\mathcal{U}$ of $M$ into $\aleph_1$ meager sets where for each ${\mathcal U}'\subset\mathcal U$,
$\bigcup {\mathcal U}'$ has the Baire property in $M$?
\end{issue}

\stepcounter{issue} 

\begin{issue}
Does there exist (in ZFC) a set of reals $X$ of cardinality $\fd$ such that all
finite powers of $X$ have Menger's property $\sfin(\cO,\cO)$?
\end{issue}

\begin{issue}
Can a Borel non-$\sigma$-compact group be generated by a Hurewicz subspace?
\end{issue}

\begin{issue}[MA]
Is there an uncountable $X\sbst\R$ satisfying $\sone(\BO,\BG)$?
\end{issue}

\begin{issue}[CH]
Is there a totally imperfect $X$ satisfying $\ufin(\cO,\Gamma)$
that can be mapped continuously onto $\Cantor$?
\end{issue}

\begin{issue}[CH]
Is there a Hurewicz $X$ such that $X^2$ is Menger but not Hurewicz?
\end{issue}

\begin{issue}
Does the Pytkeev property of $C_p(X)$ imply that $X$ has Menger's property?
\end{issue}

\begin{issue}
Does every hereditarily Hurewicz space satisfy $\sone(\BG,\BG)$?
\end{issue}

\begin{issue}[CH]
Is there a Rothberger-bounded $G\le\Bgp$ such that $G^2$ is not Menger-bounded?
\end{issue}

\begin{issue}
Let $\cW$ be the van der Waerden ideal.
Are $\cW$-ultrafilters closed under products?
\end{issue}

\begin{issue}
Is the $\delta$-property equivalent to the $\gamma$-property $\binom{\Omega}{\Gamma}$?
\end{issue}

\stepcounter{issue}

\stepcounter{issue}

\general\end{document}}

\newcommand{\Cantor}{{\{0,1\}^\N}}

\newcommand{\fb}{\mathfrak{b}}

\newcommand{\fc}{\mathfrak{c}}
\newcommand{\fd}{\mathfrak{d}}
\newcommand{\fp}{\mathfrak{p}}

\newcommand{\NON}{{\mathsf   {NON}}}
\newcommand{\COF}{{\mathsf   {COF}}}

\newcommand{\M}{\mathcal{M}}

\newcommand{\cov}{\mathsf{cov}}

\newcommand{\cf}{\mathsf{cf}}

\newcommand{\R}{\mathbb{R}}

\newcommand{\fo}{\mathfrak{od}}

\renewcommand{\split}{\mathsf{Split}}
\newcommand{\bq}{\begin{quote}}
\newcommand{\eq}{\end{quote}}
\newcommand{\cO}{\mathcal{O}}
\newcommand{\B}{\mathcal{B}}
\newcommand{\BG}{\B_\Gamma}

\newcommand{\BO}{\B_\Omega}

\newcommand{\sone}{\mathsf{S}_1}    \newcommand{\sfin}{\mathsf{S}_\mathrm{fin}}

\newcommand{\ufin}{\mathsf{U}_\mathrm{fin}}

\newcommand{\nin}{\not\in}

\newcommand{\cW}{\mathcal{W}}

\newcommand{\N}{\mathbb{N}}
\newcommand{\Z}{\mathbb{Z}}

\newcommand{\sbst}{\subseteq}
\newcommand{\by}[2]{\par\hfill\emph{#1}, #2}

\newcommand{\Tau}{\mathrm{T}}
\newcommand{\CE}{\textsc{CE}}

\newcommand{\be}{\begin{enumerate}}
\newcommand{\ee}{\end{enumerate}}
\newcommand{\bi}{\begin{itemize}}
\newcommand{\ei}{\end{itemize}}
\newcommand{\itm}{\item}

\newcommand{\general}{\small\vfill\par\noindent\hrulefill\par
\noindent\textbf{Previous issues.} The previous issues of this
bulletin are available online at\\
\texttt{http://front.math.ucdavis.edu/search?\&t=\%22SPM+Bulletin\%22}
\\[0.1cm]
\textbf{Contributions.} Announcements, discussions, and open problems should be emailed
to \texttt{tsaban@math.biu.ac.il}\\[0.1cm]
\textbf{Subscription.}
To receive this bulletin (free) to your e-mailbox, e-mail us.
}

\newcommand{\arXiv}[5]{\subsection{#2}{#4}\par\hfill{\arx{#1}}\par\hfill\emph{#3}}

\newcommand{\AMS}[3]{\subsection{#1}~\par\hfill{\texttt{#3}}\par\hfill\emph{#2}}

\newcommand{\arx}[1]{\texttt{http://arxiv.org/abs/#1}}
\newcommand{\url}[1]{\bq\texttt{#1}\eq}
\newcommand{\online}[1]{The paper is available online at \url{#1}}

\title[$\mathcal{SPM}$ Bulletin \textbf{\issuenumber} (\issuemonth{} \issueyear)]{%
$\mathcal{SPM}$ Bulletin\\[0.5cm]
Issue number \issuenumber: \issuemonth{} \issueyear{} \CE{}}

\begin{document}
\maketitle

\tableofcontents

\section{Editor's note}

This festive issue concludes the civilian year 2008 with details on a special issue of
Topology and its Applications dedicated to SPM, and with a quite large list
of research announcements.

Have a fruitful 2009,

\medskip

\by{Boaz Tsaban}{tsaban@math.biu.ac.il}

\hfill \texttt{http://www.cs.biu.ac.il/\~{}tsaban}

\section{Proceedings of the Third Workshop on Coverings, Selections and Games in Topology}

Special Issue 156 of \emph{Topology and its Applications}. Guest editor: Ljubi\v{s}a D.R. Ko\v{c}inac.

From the preface: ``The Third Workshop on Coverings, Selections and Games in Topology was held in Vrnja\v{c}ka Banja,
Serbia, from April 25 to April 29, 2007, and organized by Faculty of Sciences and Mathematics, University of Ni\v{s},
and Technical Faculty in \v{C}a\v{c}ak, University of Kragujevac.

``The previous two workshops under the same name were held in Lecce, Italy (June 27--29, 2002 and December 19--22,
2005). The main theme of this workshop was Selection Principles Theory and variety of its relations with other fields of
Topology and Mathematics (game theory, Ramsey theory, set theory, combinatorics, function spaces, hyperspaces, uniform
spaces, topological algebras, analysis, dimension theory, lattice theory, Boolean algebras).
This special issue of Topology and its Applications contains 19 papers presented at the meeting and evaluated following
the usual editorial procedure of the journal.''

Contents and links to abstracts and full texts:

\be
\itm Ljubi\v{s}a D.R. Ko\v{c}inac,
\emph{Preface}, p.~1.\\
\texttt{http://dx.doi.org/10.1016/j.topol.2008.05.022}
\itm Liljana Babinkostova,
\emph{Selective screenability in topological groups}, pp.~2--9.\\
\texttt{http://dx.doi.org/10.1016/j.topol.2008.02.014}
\itm Taras Banakh, Lyubomyr Zdomskyy,
\emph{Separation properties between the $\sigma$-compactness and Hurewicz property}, pp.~10--15.\\
\texttt{http://dx.doi.org/10.1016/j.topol.2007.12.017}
\itm Pavle V.M. Blagojevic, Aleksandra S. Dimitrijevic Blagojevic, John McCleary,
\emph{Equilateral triangles on a Jordan curve and a generalization of a theorem of Dold}, pp.~16--23.\\
\texttt{http://dx.doi.org/10.1016/j.topol.2008.04.008}
\itm Lev Bukovsk\'y,
\emph{On wQN$_*$ and wQN$^*$ spaces}, pp.~24--27.\\
\texttt{http://dx.doi.org/10.1016/j.topol.2007.10.010}
\itm Giuseppe Di Maio, Ljubi\v{s}a D.R. Ko\v{c}inac,
\emph{Statistical convergence in topology}, pp.~28--45.\\
\texttt{http://dx.doi.org/10.1016/j.topol.2008.01.015}
\itm Dragan  Djurcic, Ljubi\v{s}a D.R. Kocinac, Mali\v{s}a R. \v{Z}i\v{z}ovi\'c,
\emph{Classes of sequences of real numbers, games and selection properties}, pp.~46--55.\\
\texttt{http://dx.doi.org/10.1016/j.topol.2008.02.013}
\itm Jakub Duda, Boaz Tsaban,
\emph{Null sets and games in Banach spaces}, pp.~56--60.\\
\texttt{http://dx.doi.org/10.1016/j.topol.2008.04.009}
\itm Vitaly V. Fedorchuk, Evgenij V. Osipov
\emph{Certain classes of weakly infinite-dimensional spaces and topological games}, pp.~61--69.\\
\texttt{http://dx.doi.org/10.1016/j.topol.2007.11.007}
\itm Dimitris N. Georgiou, Stavros D. Iliadis,
\emph{On the greatest splitting topology}, pp.~70--75.\\
\texttt{http://dx.doi.org/10.1016/j.topol.2007.11.008}
\itm Stavros D. Iliadis,
\emph{Universal elements in some classes of mappings and classes of $G$-spaces}, pp.~76--82.\\
\texttt{http://dx.doi.org/10.1016/j.topol.2008.04.010}
\itm Maria Joita,
\emph{On frames in Hilbert modules over pro-$C^*$-algebras}, pp.~83--92.\\
\texttt{http://dx.doi.org/10.1016/j.topol.2007.12.015}
\itm Marion Scheepers,
\emph{Rothberger's property in all finite powers}, pp.~93--103.\\
\texttt{http://dx.doi.org/10.1016/j.topol.2007.10.011}
\itm Du\v{s}an Repov\v{s}, Boaz Tsaban, Lyubomyr Zdomskyy,
\emph{Continuous selections and $\sigma$-spaces}, pp.~104--109.\\
\texttt{http://dx.doi.org/10.1016/j.topol.2008.03.025}
\itm Andrzej Kucharski, Szymon Plewik,
\emph{Inverse systems and $I$-favorable spaces}, pp.~110--116.\\
\texttt{http://dx.doi.org/10.1016/j.topol.2007.12.016}
\itm Masami Sakai,
\emph{Function spaces with a countable cs$^*$-network at a point}, pp.~117--123.\\
\texttt{http://dx.doi.org/10.1016/j.topol.2007.10.012}
\itm Mila Mr\v{s}evi\'c, Milena Jeli\'c,
\emph{Selection principles in hyperspaces with generalized Vietoris topologies}, pp.~124--129.\\
\texttt{http://dx.doi.org/10.1016/j.topol.2008.01.016}
\itm Heike Mildenberger,
\emph{Cardinal characteristics for Menger-bounded subgroups}, pp.~130--137.\\
\texttt{http://dx.doi.org/10.1016/j.topol.2008.04.011}
\itm Tomasz Weiss,
\emph{A note on unbounded strongly measure zero subgroups of the Baer–Specker group}, pp.~138--141.\\
\texttt{http://dx.doi.org/10.1016/j.topol.2008.03.026}
\itm Sophia Zafiridou,
\emph{Dendrites with a countable closure of the set of end points}, pp.~142--149.\\
\texttt{http://dx.doi.org/10.1016/j.topol.2008.04.012}
\ee

\section{Research announcements}

\arXiv{0808.3763}
{The commutant of $L(H)$ in its ultrapower may or may not be trivial}
{Ilijas Farah, N. Christopher Phillips, Juris Stepr\=ans}
{Kirchberg asked in 2004 whether the commutant of $L(H)$ in its (norm)
ultrapower is trivial. Assuming the Continnuum Hypothesis, we
prove that the answer depends on the choice of the ultrafilter.}

\arXiv{0809.2267}
{Partitions of trees and ACA'}
{Bernard A. Anderson, Jeffry L. Hirst}
{We show that a version of Ramsey's theorem for trees for arbitrary exponents
is equivalent to the subsystem ACA' of reverse mathematics.}
{http://arxiv.org/abs/0809.2267}

\arXiv{0809.5080}
{Products of straight spaces}
{Alessandro Berarducci, Dikran Dikranjan, Jan Pelant}
{A metric space $X$ is straight if for each finite
cover of $X$ by closed sets, and for each real valued function $f$ on $X$, if
$f$ is uniformly continuous on each set of the cover, then $f$ is uniformly
continuous on the whole of $X$. A locally connected space is straight iff
it is uniformly locally connected (ULC). It is easily seen that
 ULC spaces are stable under finite products. On the
other hand the product of two straight spaces is not necessarily straight.  We
prove that the product $X\times Y$ of two metric spaces is straight if and only
if both $X$ and $Y$ are straight and one of the following conditions holds:
\begin{itemize}
\item[(a)] both $X$ and $Y$ are precompact;
\item[(b)] both $X$ and $Y$ are locally connected;
\item[(c)] one of the spaces is both precompact and locally connected.
\end{itemize}
   \par In particular, when $X$ satisfies (c), the product $X\times Z$ is straight {\em for every} straight space $Z$.
\par Finally, we characterize when  infinite products of metric  spaces are ULC and  we completely solve the problem of straightness of infinite products of ULC spaces.
}

\arXiv{0810.3030}
{Each second countable abelian group is a subgroup of a second countable divisible group}
{T. Banakh, L. Zdomskyy}
{It is shown that each pseudonorm defined on a subgroup $H$ of an abelian
group $G$ can be extended to a pseudonorm on $G$ such that the densities of the
obtained pseudometrizable topological groups coincide. We derive from this that
any Hausdorff $\omega$-bounded group topology on $H$ can be extended to a
Hausdorff $\omega$-bounded group topology on $G$. In its turn this result
implies that each separable metrizable abelian group $H$ is a subgroup of a
separable metrizable divisible group $G$. This result essentially relies on the
Axiom of Choice and is not true under the Axiom of Determinacy (which
contradicts to the Axiom of Choice but implies the Countable Axiom of Choice).
}

\arXiv{0810.1391}
{A dichotomy for Borel functions}
{Marcin Sabok}
{The dichotomy discovered by Solecki states that any Baire class
1 function is either $\sigma$-continuous or "includes" the Pawlikowski function
$P$. The aim of this paper is to give an argument which is simpler than the
original proof of Solecki and gives a stronger statement: a dichotomy for all
Borel functions.}

\arXiv{0810.5587}
{More Results on Regular Ultrafilters in ZFC}
{Paolo Lipparini}
{We prove, in ZFC alone, some new results on regularity and decomposability of
ultrafilters.
 We also list some problems, and furnish applications to topological spaces
and to extended logics.}

\arXiv{0811.0914}
{Minimal pseudocompact group topologies on free abelian groups}
{Dikran Dikranjan, Anna Giordano Bruno, Dmitri Shakhmatov}
{A Hausdorff topological group $G$ is minimal if every continuous isomorphism $f:G\to H$ between $G$ and a Hausdorff topological group $H$ is
open. Significantly strengthening a 1981 result of Stoyanov we prove the following theorem: For every infinite minimal group $G$ there exists a sequence $\{\sigma_n:n\in\N\}$ of cardinals such that
$$
w(G)=\sup\{\sigma_n:n\in\N\}\ \ \mbox{ and } \ \
\sup\{2^{\sigma_n}:n\in\N\}\leq |G| \leq 2^{w(G)},
$$
where $w(G)$ is the weight of $G$. If $G$ is an infinite minimal abelian group, then either $|G|=2^\sigma$ for some cardinal $\sigma$, or $w(G)=\min\{\sigma:  |G|\le 2^\sigma\}$. Moreover, we show that the equality $|G|=2^{w(G)}$ holds whenever $\cf(w(G))> \omega$.

For a cardinal $\kappa$ we denote by $F_\kappa$ the free abelian group with $\kappa$ many generators. If $F_\kappa$ admits a pseudocompact group topology, then $\kappa\ge\fc$, where $\fc$ is the cardinality of the continuum. We show that the existence of a minimal pseudocompact group topology  on $F_\fc$ is equivalent to the Lusin's Hypothesis $2^{\omega_1}=\mathfrak c$.  For $\kappa>\mathfrak c$, we prove that  $F_\kappa$  admits a (zero-dimensional) minimal pseudocompact group topology if and only if $F_\kappa$ has both a minimal group topology and a pseudocompact group topology.
If $\kappa>\fc$, then $F_\kappa$  admits a connected minimal pseudocompact group topology of weight $\sigma$ if and only if $\kappa=2^\sigma$. Finally, we establish that no infinite torsion-free abelian group can be equipped with a locally connected minimal group topology.}

\arXiv{0811.3661}
{Weakly infinite dimensional subsets of $\R^\N$}
{Liljana Babinkostova and Marion Scheepers}
{The Continuum Hypothesis implies an Erd\"os-Sierpi\'nski like duality between
the ideal of first category subsets of $\R^{\N}$, and the ideal of
countable dimensional subsets of $\R^{\N}$. The algebraic sum of a
Hurewicz subset - a dimension theoretic analogue of Sierpinski sets and Lusin
sets - of $\R^{\N}$ with any strongly countable dimensional subset
of $\R^{\N}$ has first category.
}

\arXiv{0811.3945}
{Selections, Extensions and Collectionwise Normality}
{Valentin Gutev and Narcisse Roland Loufouma Makala}
{We demonstrate that the classical Michael's selection theorem for l.s.c.
mappings with a collectionwise normal domain can be reduced only to
compact-valued mappings modulo the Dowker's extension theorem for such spaces.
The technique developed to achieve this result is applied to construct
selections for set-valued mappings whose point images are in completely
metrizable absolute retracts.}

\arXiv{0811.4272}
{Openly factorizable spaces and compact extensions of topological semigroups}
{Taras Banakh, Svetlana Dimitrova}
{We prove that the semigroup operation of a topological semigroup $S$ extends
to a continuous semigroup operation on its the Stone-\v{C}ech compactification
$\beta S$ provided $S$ is a pseudocompact openly factorizable space, which
means that each map $f:S\to Y$ to a second countable space $Y$ can be written
as the composition $f=g\circ p$ of an open map $p:X\to Z$ onto a second
countable space $Z$ and a map $g:Z\to Y$. We present a spectral
characterization of openly factorizable spaces and establish some properties of
such spaces.}

\arXiv{0811.4276}
{Embedding the bicyclic semigroup into countably compact topological semigroups}
{Taras Banakh, Svetlana Dimitrova, Oleg Gutik}
{We study algebraic and topological properties of topological semigroups
containing a copy of the bicyclic semigroup $C(p,q)$. We prove that each
topological semigroup $S$ with pseudocompact square contains no dense copy of
$C(p,q)$. On the other hand, we construct a consistent example of a Tychonov
countably compact semigroup containing a copy of $C(p,q)$.}

\arXiv{0812.0489}
{Building suitable sets for locally compact groups by means of continuous selections}
{Dmitri Shakhmatov}
{If a discrete subset $S$ of a topological group $G$ with the identity $1$ generates
a dense subgroup of $G$ and $S \cup {1}$ is closed in $G$, then $S$ is called a
suitable set for $G$. We apply Michael's selection theorem to offer a direct,
self-contained, purely topological proof of the result of Hofmann and Morris on
the existence of suitable sets in locally compact groups. Our approach uses
only elementary facts from (topological) group theory.}

\AMS{The near coherence of filters principle does not imply the filter dichotomy principle}
{Heike Mildenberger and Saharon Shelah}
{http://www.ams.org/journal-getitem?pii=S0002-9947-08-04806-X}

\ed